\documentclass[11pt,eqno]{article}
\usepackage{amssymb}
\begin{document}

\newcommand{\defi}{\stackrel{\Delta}{=}}
\newcommand{\qed}{\hphantom{.}\hfill $\Box$\medbreak}
\newcommand{\A}{{\cal A}}
\newcommand{\B}{{\cal B}}
\newcommand{\U}{{\cal U}}
\newcommand{\G}{{\cal G}}
\newcommand{\cZ}{{\cal Z}}
\newcommand{\proof}{\noindent{\bf Proof \ }}
\newcommand\one{\hbox{1\kern-2.4pt l }}
\newcommand{\Item}{\refstepcounter{Ictr}\item[(\theIctr)]}
\newcommand{\QQ}{\hphantom{MMMMMMM}}

\newtheorem{Theorem}{Theorem}[section]
\newtheorem{Lemma}[Theorem]{Lemma}
\newtheorem{Corollary}[Theorem]{Corollary}
\newtheorem{Remark}[Theorem]{Remark}
\newtheorem{Proposition}[Theorem]{Proposition}
\newtheorem{Definition}[Theorem]{Definition}
\newtheorem{Construction}[Theorem]{Construction}
\newtheorem{Exa}[Theorem]{Example}

\newcounter{Ictr}

\renewcommand{\theequation}{
\arabic{equation}}
\renewcommand{\thefootnote}{\fnsymbol{footnote}}

\def\C{\mathcal{C}}

\def\V{\mathcal{V}}

\def\I{\mathcal{I}}

\def\Y{\mathcal{Y}}

\def\X{\mathcal{X}}

\def\J{\mathcal{J}}

\def\Q{\mathcal{Q}}

\def\W{\mathcal{W}}

\def\S{\mathcal{S}}

\def\T{\mathcal{T}}

\def\L{\mathcal{L}}

\def\M{\mathcal{M}}

\def\N{\mathcal{N}}
\def\R{\mathbb{R}}
\def\H{\mathbb{H}}


\begin{center}
\topskip3cm

 {\bf Homogeneous
Cone Complementarity Problems and  $P$ Properties
 \footnote{The work was supported in part by the National Natural
 Science Foundation of China (10831006) and the National Basic Research Program of China (2010CB732501), and a Discovery Grant from NSERC.} }

\vskip3mm
 \renewcommand{\thefootnote}{\fnsymbol{footnote}}
 Lingchen Kong\footnote{~Lingchen Kong,
Department of Applied Mathematics, Beijing Jiaotong
University, Beijing 100044, P. R. China (e-mail:
konglchen@126.com)}, Levent Tun\c{c}el\footnote{~Levent
Tun\c{c}el, Department of Combinatorics and Optimization, Faculty
of Mathematics, University of Waterloo, Waterloo, Ontario N2L 3G1,
Canada (e-mail: ltuncel@math.uwaterloo.ca)}, and Naihua
Xiu\footnote{~Naihua Xiu, Department of Applied Mathematics,
Beijing Jiaotong University, Beijing 100044, P. R. China (e-mail: nhxiu@center.njtu.edu.cn)}\\
\vskip4mm
(April~9,~2009; Revised November 8, 2010)\\
\end{center}
\vskip12pt

\begin{abstract} We consider existence and uniqueness
properties of a solution to homogeneous cone complementarity problem
(HCCP). Employing the $T$-algebraic characterization of homogeneous
cones, we generalize the $P, P_0, R_0$ properties for a nonlinear
function associated with the standard nonlinear complementarity
problem to the setting of HCCP.  We prove that if a continuous
function has either the order-$P_0$ and $R_0$, or the $P_0$ and
$R_0$ properties then all the associated HCCPs have solutions. In
particular, if a continuous monotone function has the trace-$P$
property then the associated HCCP has a unique solution (if any); if
it has the uniform-trace-$P$ property then the associated HCCP has
the global uniqueness (of the solution) property (GUS). Moreover, we
establish a global error bound for the HCCP with the
uniform-trace-$P$ property under some conditions
for homogeneous cone linear complementarity problem.

\noindent

\vskip10pt

 \vskip10pt


\vskip10pt

\end{abstract}

\section{Introduction}
In this paper, we are interested in the \emph{homogeneous cone
complementarity problem} (HCCP$(F,q)$ for short) which is to find a
vector $x\in K$ such that
\begin{eqnarray}\label{100}~x \in \textit{K}, ~y \in {K^\ast},~\langle
x,y\rangle =0, ~~y=F(x)+q, \end{eqnarray} where  $K$ is a
homogeneous cone (the automorphism group of the cone acts
transitively on the interior of the cone, see Section 2 for the details)
in a finite-dimensional inner product space $\H$ over $\R$ with its
dual $K^\ast$ given by $K^\ast:=\{y\in \H: \langle x,y\rangle\geq 0,
\forall x\in K\},$  $F:\H\rightarrow \H$ is a continuous function
and $q\in\H$. If $F(x)=L(x)$ is linear, we call problem (\ref{100})
the \emph{homogeneous cone linear complementarity problem}
(HCLCP$(L,q)$). When $K$ is a symmetric cone in a Euclidean Jordan
algebra, it is the \emph{symmetric cone complementarity problem}
(SCCP), which includes the so-called \emph{nonlinear complementarity
problem} (NCP, where $\H=\R^n$, the space of $n$-dimensional real
column vectors,  and $K=\R^n_+$, the nonnegative orthant) and
\emph{semidefinite complementarity problem} (SDCP, where
$\H={\mathbb S}^n$, the space of $n \times n$ real symmetric
matrices, and $K={\mathbb S}^n_+$, the cone of positive semidefinite
symmetric matrices) as special cases.

\vskip1.5mm In the NCP context, a continuous function
$f:\R^n\rightarrow\R^n$ is said to be a \emph{$P$ function (has
the $P$ property)} if the following implication holds
$$(x-y)\circ [f(x)-f(y)]\leq 0~~~\Rightarrow~~~x=y,$$
where $``\circ"$ denotes the Hadamard (componentwise) product and
$z\leq 0$ means that all components of $z$ are nonpositive. There
are many applications of $P$ functions in engineering, economics,
management science, as well as several
other fields, see, e.g., \cite{FP03, Sch04}.
We say $f$ is a $P_0$ function if $f+\varepsilon I$ is a $P$
function for any $\varepsilon>0$ where $I$ is the identity
transformation. This is a generalization of $P$-matrices. It is
known that if $f$ is a $P_0$ function and satisfies the so-called
$R_0$ condition then the NCP$(f,q)$ has a solution for every $q\in
\R^n$. (Here, we only consider the $R_0$ condition as described in
Definition \ref{d01} below.) In the setting of a Cartesian product
of sets in $\R^n$, Facchinei and Pang \cite{FP03} proposed
the notions of $P$ and $P_0$
functions and studied some of their properties. In the setting of
Euclidean Jordan algebras, Gowda, Sznajder and Tao \cite{GST04}
studied some $P$ and $P_0$ properties for linear transformations;
Tao and Gowda \cite{TG05} introduced $P$ and $P_0$ functions and
established the existence result for SCCP. Moreover,  Gowda and
Sznajder \cite{GS06} studied the automorphism invariance of $P$ and
\emph{globally uniquely
solvability} (GUS) properties for linear transformations
on Euclidean Jordan algebras.  
a necessary condition for $F$
Symmetric cones are \emph{homogeneous}  and \textit{self-dual},
see \cite{FK94, K99}. A natural next step in this
line of generalizations is to
drop the requirement that $K$ is self-dual. While there is a
finite number of non-isomorphic symmetric cones of each dimension,
the number is uncountable for homogeneous cones when the dimension
$n>11$, see \cite{Vin63}. There has been some
increase in the interest and
activity in the area of homogeneous cones and optimization
problems over homogeneous cones, see, e.g.,
\cite{AW04,Chua03,Chua08,Do79,D279,Fay02,Gu96,GT98,Is00,Ost98,Ro66,TX01,TT04,Vin63,Vin65}.
These papers deal with either certain theoretical properties of
homogeneous cones, primal-dual interior-point methods for linear
programming over homogeneous cones or their applications. In this
paper, we work on the $P$ properties for nonlinear transformations
in the setting of HCCP. The aim of our work is to establish the
existence and uniqueness results of a solution to HCCP.

\vskip1.5mm  With the help of the $T$-algebraic characterization of
homogeneous cones, we first study the \emph{metric projection} onto
homogeneous cone $K$ and its properties related to HCCP. Based on
them, we introduce $P$, order-$P$, trace-$P$, uniform-trace-$P$,
trace-$P_0$,  order-$P_0$, $P_0$ and $R_0$ properties for a function
$F:\H\rightarrow\H$ in the setting of HCCP. Then, we show that if
$F$ has either the order-$P_0$ and $R_0$, or  the $P_0$ and $R_0$
properties then the HCCP$(F,q)$ has a solution for every $q\in \H$
by applying degree theory. We give some sufficient conditions
under which the HCCP has GUS property. For instance, if $F$ is
monotone and has the trace-$P$ property then the associated
HCCP$(F,q)$ has a unique solution (if not empty); if $F$ is monotone
and has the uniform-trace-$P$ property then the associated
HCCP$(F,q)$ has GUS property. Moreover, we establish a global error
bound for the HCCP with $F$ having the uniform-trace-$P$ property.

\vskip1.5mm This paper is organized as follows. In Section 2, we
briefly review some basic concepts and results on $T$-algebras, and
describe some fundamental results on metric projection onto
homogenous cones. In Section 3, we introduce various $P$
properties and show our existence result for HCCP. In Section 4, we
study the GUS property and give an error bound for HCCP.  In Section 5, we include some concluding
remarks. While most of our focus is on Homogeneous cones, many of
our results apply more generally (in the setting of arbitrary convex
cones).

\section{Preliminaries}

We first briefly review some basic concepts and results on
homogeneous cones and $T$-algebras from \cite{Chua08,Vin63,Vin65},
and then provide some fundamental results on metric projection onto
homogenous cones.

\subsection{Homogeneous cones and $T$-algebras}

\begin{Definition}\label{d1} A closed, convex cone $K$ with nonempty
interior is \emph{homogeneous} if the group of automorphisms of
$K$ acts transitively on the interior of $K$.
\end{Definition}

Note that a cone $K$ is homogeneous then so is its dual $K^\ast$.
Vinberg \cite{Vin63} introduced a constructive way to build
homogeneous cones by employing the so-called \emph{$T$-algebra}
which connects homogeneous cones to abstract matrices whose elements
are vectors. We first  review the following concept of \emph{matrix
algebra}.

\begin{Definition}\label{d2} A \emph{matrix algebra} $\A$ is a
bi-graded algebra $\bigoplus_{i,j=1}^r \A_{ij}$ over the reals
with  a bilinear product of  $a_{ij}\in\A_{ij}$ and
$a_{kl}\in\A_{kl}$ ($1\leq i, j,k, l \leq r$) satisfying
\begin{eqnarray}\label{101}
a_{ij}a_{kl}\in\left \{\begin{array}{cc} \A_{il}  &
\textrm{if}~~j=k,\cr\noalign{\vskip2truemm}
~0  & \textrm{if}~~j\neq k, \\
\end{array}\right.
\nonumber\end{eqnarray} where $\A_{ij}$ is a $n_{ij}$-dimensional
vector space. The positive integer $r$ is called the \emph{rank}
of the matrix algebra $\A$.
\end{Definition}

Every element $a\in\A$ is a generalized matrix with its component
in $\A_{ij}$ being an $n_{ij}$-dimensional generalized element of
the matrix $a_{ij}$, i.e., $a_{ij}$ is the projection of $a$ onto
$\A_{ij}$. The multiplication of two elements $a,b\in\A$ is
analogous to the multiplication of matrices,
\begin{eqnarray}\label{102} (ab)_{ij} =\sum^r_{k=1}
a_{ik}b_{kj}.\nonumber\end{eqnarray} Assume that for every $i$,
$\A_{ii}$ is isomorphic to $\R$, and let $\rho_i$ be the
isomorphism and let $e_i$ denote the representation of the \emph{unit
element of $\A_{ii}$ in $\A$}. We define the \textit{trace} of an
element $a$ as
\begin{eqnarray} \label{103} \textup{Tr} (a) = \sum^r_{i=1} \rho_i
(a_{ii}).\nonumber\end{eqnarray}

The following notion generalizes the classical (conjugate)
transpose. An \emph{involution} $\ast$ of the matrix algebra $\A$ of
rank $r$ is a linear automorphism on $\A$ that satisfies

(i) $(a^\ast)^\ast=a$ for all $a\in\A$ (involutory).

(ii) $(ab)^\ast=b^\ast a^\ast$ for all $a,b\in\A$
(anti-homomorphic).

(iii) $(a^\ast)_{ij} = (a_{ji})^\ast$  for $1\leq i, j  \leq r$.

(iv) $\A^\ast_{ij}= \A_{ji}$ for $1\leq i, j  \leq r$.

We are ready to state the following definition of $T$-algebra, which
was originally introduced by Vinberg \cite{Vin63}.

\begin{Definition}\label{d3} A \emph{$T$-algebra} of rank $r$ is a matrix algebra $\A$ of rank $r$ with
involution $\ast$ satisfying the following axioms:

I. For each $1\leq i\leq r$, the subalgebra $\A_{ii}$ is
isomorphic to the reals.

II. For each $a\in\A$ and each $1\leq i, j\leq r$, $$ a_{ji}e_i =
a_{ji}~~{\rm and}~~ e_ia_{ij} = a_{ij}.$$

III. For each $a,b\in\A$, \textup{Tr}$(ab)= $\textup{Tr}$(ba)$.

IV. For each $a, b, c\in\A$, \textup{Tr}$((ab)c)=
$\textup{Tr}$(a(bc))$.

V. For each $a\in\A$,  \textup{Tr}$(a^\ast a)\geq 0$, with
equality if and only if $a=0$.

VI. For each $a, b, c\in\A$ and each $i,j,k,l\in\{1,2,\cdots,r\}$
with $i\leq j\leq k\leq l$,
$$a_{ij}(b_{jk}c_{kl})=(a_{ij}b_{jk})c_{kl}.$$

VII. For each $a,b\in\A$ and each $i,j,k,l\in\{1,2,\cdots,r\}$
with $i\leq j\leq k$ and $l\leq k$, $$a_{ij}(b_{jk}b^\ast_{lk})=
(a_{ij}b_{jk})b^\ast_{lk}.$$
\end{Definition}
Thus, the \textit{cone associated with a $T$-algebra} $\A$ of rank
$r$, denoted by ${\rm int}(K(\A))$, is given by
\begin{eqnarray}\label{104} {\rm int}(K(\A)):=\{tt^\ast: t\in\A,  t_{ij}=0,
~\forall 1\leq j<i\leq r,~ {\rm and}~ t_{ii}>0, ~\forall 1\leq i\leq
r\}.\nonumber\end{eqnarray} It is easy to see from Axiom II that
$e:= \sum_{i=1}^r e_i$ is the unit element in $\A$,
i.e., $ea=ae$ for all $a \in \A$.  It is necessary to note that
multiplication in a $T$-algebra is neither commutative nor
associative. Define the subalgebra of upper triangular elements of
$\A$ and the subspace of ``Hermitian" elements, respectively, as
\begin{eqnarray}\label{105}\T:= \bigoplus_{i\leq j, i,j=1}^r
\A_{ij},~~~~~~~{\cal H}: = \{a\in\A : a = a^\ast\}.
\nonumber\end{eqnarray} Clearly,  Axiom VI is equivalent to $t(uw)
= (tu)w$ for all $t,u,w\in\T$. Taking involution, we get another
equivalent statement: $t(uw) = (tu)w$ for all $t,u,w\in\T^\ast$.
Similarly, Axiom VII is equivalent to $t(uu^\ast) = (tu)u^\ast$
for all $t, u\in\T$; or  equivalently, $(u^\ast u)t = u^\ast(ut)$
for all $t,u\in\T^\ast$.

Let $\T_\ast$ (resp., $\T_+$ and $\T_{++}$) denote the set of
elements of $\T$ with nonzero (resp., nonnegative and positive)
diagonal components. It is easy to see that ${\rm
int}(K(\A))=\{tt^\ast: t\in\T_{++}\}$ and ${\cal H}$ is the linear
span of ${\rm int}(K(\A))$.

We recall below a fundamental characterization of homogeneous cones
established by Vinberg \cite{Vin63}.

\begin{Theorem}\label{t1} ($T$-algebraic representation of homogeneous cones) A cone $K$ is
homogeneous if and only if  ${\rm int}(K)$ is isomorphic to the cone
${\rm int}(K(\A))$ associated with some $T$-algebra $\A$. Moreover,
given ${\rm int}(K(\A))$, the representation of an element from
${\rm int}(K)$ in the form $tt^\ast$ is unique. Finally, the dual
cone ${\rm int}(K^\ast)$ can be represented as $\{t^\ast
t:t\in\T_{++}\}$.
\end{Theorem}

\noindent\textbf{Remark 2.1} This theorem is analogous to the
representation of a symmetric cone as the set of squares over a
Jordan algebra. However, for an arbitrary
$T$-algebra and its associated cone it
is not true that for a given $a\in\A$ we have $aa^\ast\in {\rm
int}(K(\A))$.
In general, the closures  of
the cone ${\rm int}(K(\A))$ and its dual are specified by
\begin{eqnarray}\label{106} K(\A)=\{tt^\ast:
t\in\T_{+}\}, ~~K^\ast(\A)=\{t^\ast t: t\in\T_{+}\}.\end{eqnarray}
Thus, $K(\A)$ is a closed convex cone in the inner product space
$({\cal H}, \langle \cdot,\cdot\rangle)$, where $\langle
\cdot,\cdot\rangle$ is given by $\langle a,
b\rangle := \textup{Tr}(a^\ast b)=\textup{Tr}(ab)$ for all $a,b\in
{\H}$. It is worth noting that for all $a,b,c\in\A$,
\begin{eqnarray}\label{t0t}\langle ab,c\rangle
=\langle b^\ast a^\ast,c^\ast\rangle
=\langle a, cb^\ast\rangle=\langle b,a^\ast c\rangle.\end{eqnarray}
We define the \emph{norm} induced by the inner product as
$\|a\|:=\sqrt{\langle a,a\rangle}$. In what follows, we simply
write $K$ and $\H$ for $K(\A)$ and $({\cal H}, \langle
\cdot,\cdot\rangle)$, respectively.

\vskip1.5mm\noindent\textbf{Example 2.1} Consider the following
five-dimensional closed convex cone with nonempty interior (Vinberg
\cite{Vin63}):
$$K:=\left\{x\in\R^5:\left(
                    \begin{array}{ccc}
                      x_1 & x_2 & x_4 \\
                      x_2 & x_3 & 0 \\
                     x_4 & 0 & x_5 \\
                    \end{array}
                  \right)\in{\mathbb S}_+^3\right\}.$$
This cone is homogeneous; but, it is not a symmetric cone since
there does not exist any inner product on $\R^5$ under which
$K=K^\ast$. Let us choose the inner product implied by the trace
inner product on ${\mathbb S}^3$; that is, $\forall x,y\in \R^5$,
$$\langle x,y\rangle:=x_1y_1+x_3y_3+x_5y_5+2x_2y_2+2x_4y_4.$$
Then,
$$K^\ast=\left\{y\in\R^5: \left(
                                    \begin{array}{cc}
                                      y_1 & y_2 \\
                                      y_2 & y_3 \\
                                    \end{array}
                                  \right)\in{\mathbb S}^2_+,\left(
                                    \begin{array}{cc}
                                      y_1 & y_4 \\
                                      y_4 & y_5 \\
                                    \end{array}
                                  \right)\in{\mathbb S}^2_+\right\}.$$
With this (natural) choice of the inner product,  $K\subseteq
K^\ast$. Moreover, it is straightforward to verify that $a^2\in K+K^\ast$
for all $a\in \H$.

\vskip1.5mm\noindent\textbf{Remark 2.2}  Note that every $z\in\H$
may be rewritten as $z=t+t^\ast$ with $t\in\T$ and $\langle
z^2,e_i\rangle\geq 0$ for all $i\in\{1,2,\cdots,r\}$. Then
$z^2=tt^\ast+t^\ast t+t^2+(t^\ast)^2$. Motivated by the fact that a
symmetric cone is the set of squares over a Jordan algebra, we
propose the following question:
\begin{eqnarray}
\label{Q1}
\begin{tabular}{|c|}\hline
\parbox{4.2in}{
\[
\mbox{for every homogeneous cone $K$, does there exist }
\]
\[
\mbox{an inner product on $\H$ such that $z^2\in K+K^\ast, \forall
z\in \H?$}
\]
}
\\ \hline
\end{tabular}
\nonumber
\end{eqnarray}

As we stated in the Introduction, we are generalizing the underlying
theorems and in many cases their existing proofs from the symmetric
cone setting to the more general, homogeneous cone setting. The
existing results for the symmetric cone setting essentially fix an
inner product under which $K=K^\ast$ and treat both $K$ and $K^\ast$
in the same finite dimensional Euclidean space. This allows
operations like $x+y$ for $x\in K, y\in K^\ast$ and the related
metric projections onto the cones $K,K^\ast$ to be treated in the same
space. Since much of the related theory is based on metric
projections, one is required to fix an inner product in our more
general setting as well. We remind the reader that the choice of the
inner product is up to the goals of the user of the theory (for
example, to recover the existing results for the special case of
symmetric cones, we would pick the inner product so that
$K=K^\ast$).

\subsection{Metric Projection}
Let $\Pi_K(x)$ denote the \emph{metric projection} of $x$ onto
$K$, i.e.,
\begin{eqnarray} \label{100-1}\Pi_K(x):= {\rm
argmin}\left\{\frac{1}{2}\|x-z\|^2: z\in
K\right\}.\nonumber\end{eqnarray} In other words, $y=\Pi_K(x)$ if
and only if $y\in K$  and
\begin{eqnarray}\label{100-1b}\|x-y\|\leq\|x-z\|, ~\forall z\in
K, \nonumber\end{eqnarray} or equivalently,  the so-called obtuse
angle property (or the Kolmogorov criterion) holds:
\begin{eqnarray}\label{100-2}\langle z-\Pi_K(x), x-\Pi_K(x)\rangle\leq 0, ~\forall z\in
K. \nonumber\end{eqnarray} It is well-known \cite{Zar71} that the
metric projector $\Pi_K$ is unique and contractive, i.e.,
\begin{eqnarray}\label{100-3}\|\Pi_K(x)-\Pi_K(y)\|\leq\|x-y\|,~\forall x,y\in\H.\nonumber\end{eqnarray}
Utilizing the Moreau decomposition,  any $x\in\H$ can be written as
\begin{eqnarray}\label{00} x= \Pi_K(x)-\Pi_{K^\ast}(-x)~~{\rm
with}~\langle\Pi_K(x),\Pi_{K^\ast}(-x)\rangle = 0.\end{eqnarray}

Based on the metric projection operator, we define the
following operations for any $x,y\in\H$,
\begin{eqnarray}\label{001} x\wedge_K y&:=&x-\Pi_K(x-y),~~~~
 x\vee_K y~:=~y+\Pi_K(x-y).\end{eqnarray} Then, by direct
calculation, we obtain
\begin{eqnarray}x\wedge_K y
=y\wedge_{K^\ast}x,~~~~
(-x)\wedge_K (-y)
=-(x\vee_{K^\ast} y).\end{eqnarray}

Summarizing the
above arguments, we have the following proposition.
\begin{Proposition} \label{t101}
Let $K$ be a closed convex cone in $\H$ with its dual $K^\ast$. Then
the following statements hold for all $x,y\in\H$.
\begin{description}
  \item [](a)~~We have  $x= \Pi_K(x)-\Pi_{K^\ast}(-x)~~{\rm
with}~\langle\Pi_K(x),\Pi_{K^\ast}(-x)\rangle = 0$. This
decomposition is unique in the sense that if $x = x_1- x_2$ with
$x_1\in K, x_2 \in K^\ast$ and $\langle x_1, x_2\rangle=0$ then
$x_1 = \Pi_K(x)$ and $x_2 = \Pi_{K^\ast}(-x)$.
  \item [](b)~~$x\wedge_K y=y\wedge_{K^\ast}x,~~x\vee_K y=y\vee_{K^\ast}x.$
   \item [](c)~~$(-x)\wedge_K(-y)=-(x\vee_{K^\ast}y)$ and ~$(-x)\wedge_{K^\ast}(-y)=-(x\vee_Ky).$
\end{description}
In particular, $\Pi_K(x)\wedge_K\Pi_{K^\ast}(-x)=0,$ and $
~\Pi_K(x)\vee_K\Pi_{K^\ast}(-x)=\Pi_K(x)+\Pi_{K^\ast}(-x).$
\end{Proposition}

Considering the characterization of homogeneous cones, from
(\ref{00}) and (\ref{106}), we obtain that any $x\in\H$ can be
expressed as
\begin{eqnarray}\label{107-0}x=uu^\ast-v^\ast v~~{\rm with}~~\langle
uu^\ast,v^\ast v\rangle=0,\nonumber\end{eqnarray} where
$u,v\in\T_{+}$. Observe that
\begin{eqnarray} \langle
uu^\ast,v^\ast v\rangle &=& \langle
(uu^\ast)v^\ast, v^\ast\rangle~~~\nonumber\\
&=& \langle
u(u^\ast v^\ast),v^\ast \rangle~~~\nonumber\\
&=& \langle
(vu)^\ast,u^\ast v^\ast \rangle \nonumber\\
&=& \langle(vu)^\ast,(vu)^\ast\rangle\nonumber\\
&=& \|vu\|^2, ~~ \nonumber\end{eqnarray}where the first equality
holds by (\ref{t0t}), the second equality holds by  Axiom VII, the
third holds by (\ref{t0t}) and the fact that $\ast$ is
anti-homomorphic. Then, Axiom V implies that $\langle uu^\ast,v^\ast
v\rangle=0$ if and only if $vu=0$. We have actually proved the
following.

\begin{Theorem}\label{0t1}
Let $K$ be a homogeneous cone in $\H$ with its dual $K^\ast$.
Then, every $x\in{\H}$ can be uniquely
expressed as
\begin{eqnarray}\label{107}x=uu^\ast-v^\ast
v~~{\rm with}~~vu=0,~~u,v\in\T_{+}.\end{eqnarray}
Moreover, we have $\Pi_{K}(x)=uu^\ast,~~ \Pi_{K^\ast}(-x)=v^\ast
v.$
\end{Theorem}

Applying Proposition \ref{t101} and the above theorem, we obtain
the following equivalent statements related to HCCP (\ref{100}).
\begin{Proposition} \label{t102-2} Let $K$ be
a homogeneous cone in ${\H}$ with its dual $K^\ast$. Then
the following statements are equivalent:
\begin{description}
  \item [](a)~~$x\wedge_K y=0;$
   \item [](b)~~$y\wedge_{K^\ast}x=0;$
   \item [](c)~~$x\in K,~y\in K^\ast,~\langle x,y\rangle=0;$
  \item [](d)~~$x\in K,~y\in K^\ast,~\langle xy, e_i\rangle
=\langle yx, e_i\rangle=0,~\forall i\in\{1,2,\cdots,r\};$
   \item [](e)~~There exist $u,v\in\T_{+}$ such that
$x=uu^\ast,~y=v^\ast v,~~vu=0;$
 \item [](f)~~$x\in K,~y\in K^\ast,~(xy)_{lj}=0,~\forall~  l,j\in\{1,2,\cdots,r\}~such~that~l\geq j.$
\end{description}
In particular,  if $xy=yx$, then (f) becomes

\vskip1mm $(f')$~~$x\in K,~y\in K^\ast,~xy=0.$

\end{Proposition}
\textbf{Proof.} By Proposition \ref{t101} and Theorem \ref{0t1},
$(a)\Leftrightarrow (b) \Leftrightarrow (c) \Leftrightarrow (e)$.
Clearly, $(f)\Rightarrow (d)\Rightarrow (c)$. Therefore, we need
only to show that $(e)\Rightarrow (f)$. Choose any $w^\ast\in
\T^{\ast}$. Applying arguments similar to those before Theorem
\ref{0t1}, we obtain
\begin{eqnarray}  \langle
xy,w^\ast\rangle&=&\langle
(uu^\ast)(v^\ast v), w^\ast\rangle \nonumber\\
&=& \langle
uu^\ast, w^\ast(v^\ast v)\rangle ~~~\nonumber\\
&=& \langle
uu^\ast, (w^\ast v^\ast) v\rangle ~~~~\nonumber\\
&=&\langle
(uu^\ast)v^\ast, w^\ast v^\ast\rangle ~~~~\nonumber\\
&=& \langle u(vu)^\ast, w^\ast
v^\ast\rangle.~~~~\nonumber\end{eqnarray} Then, the desired
conclusion of equivalence follows. $(f')$ is self-evident.{\qed}

We next address the following result which will be used to establish
the GUS-property for HCCP.
\begin{Proposition} \label{t103}
Let $K$ be a homogeneous cone in $\H$ with its dual $K^\ast$.
Then, for every $x,y\in \H$,
the following statements hold:
\begin{description}
  \item [](i)~~if $x\in K,~y\in K^\ast,$ and $x=\sum_{i=1}^rx_ie_i$, then $\langle xy, e_i\rangle\geq 0,~\forall i\in\{1,2,\cdots, r\};$
  \item [](ii)~$\langle (x\wedge_K y)(x\vee_{K}y),e_i\rangle=\langle xy, e_i\rangle,~\forall i\in\{1,2,\cdots, r\}.$
\end{description}
\end{Proposition}
\textbf{Proof.}   We first prove (i). Since $x\in K,~y\in K^\ast$,
noting that $e_i\in K\cap
K^\ast$,  we have $\langle x,e_i\rangle\geq0, \langle y,e_i\rangle\geq0$ for
all $i\in\{1,2,\cdots,r\}$. Using
$x=\sum_{i=1}^rx_ie_i$, $x_i=\langle x,e_i\rangle\geq0$.  Thus,
for every $i\in\{1,2,\cdots,r\}$,
\begin{eqnarray}  \langle
xy,e_i\rangle=\langle y, x e_i\rangle = x_i\langle y, e_i\rangle
\geq 0.~~~~\nonumber\end{eqnarray} Therefore, the desired conclusion
(i) holds.

\vskip1mm For part (ii),  direct calculation yields
\begin{eqnarray}\langle (x\wedge_K y)(x\vee_{K}y), e_i\rangle &=&\langle
[x-\Pi_K(x-y)][y+\Pi_K(x-y)], e_i\rangle\nonumber\\
&=&\langle xy,e_i\rangle+\langle x\Pi_K(x-y),e_i\rangle-\langle
\Pi_K(x-y)y,e_i\rangle-\langle \Pi_K(x-y) \Pi_K(x-y),e_i
\rangle\nonumber\\
&=&\langle xy,e_i\rangle+\langle x,\Pi_K(x-y)e_i\rangle-\langle
y,\Pi_K(x-y)e_i\rangle-\langle \Pi_K(x-y) ,\Pi_K(x-y)e_i
\rangle\nonumber\\&=&\langle xy,e_i\rangle+\langle (x-y)-\Pi_K(x-y), \Pi_K(x-y)e_i\rangle\nonumber\\
&=&\langle xy,e_i\rangle,\nonumber\end{eqnarray} where the last
equality follows from the fact $\langle (x-y)-\Pi_K(x-y),
\Pi_K(x-y)e_i\rangle=\langle \Pi_{K^\ast}(y-x)
\Pi_K(x-y),e_i\rangle=0$ by (\ref{00}).  Thus, we
proved (ii). {\qed}

Note that in the above proposition, the condition
$x=\sum_{i=1}^rx_ie_i$ is necessary, which is illustrated by the
following example.

\vskip1.5mm\noindent\textbf{Example 2.1} Let $K$ be the homogeneous
cone in $\R^5$ as in Example 2.1. Clearly, it is easy to verify that
$$x=\left(
                    \begin{array}{ccc}
                      5 & -2 & -2 \\
                      -2 & 1 & 0 \\
                     -2 & 0 & 5 \\
                    \end{array}
                  \right)\in K,~~
y=\left(
                    \begin{array}{ccc}
                      1 &  2 &  2 \\
                       2 & 4 & 0 \\
                      2 & 0 & 4 \\
                    \end{array}
                  \right)\in K^\ast,~~{\rm and}~~ x y=\left(
                    \begin{array}{ccc}
                      -3 &  2 &  2 \\
                       0 & 0 & 0 \\
                      8 & 0 & 16 \\
                    \end{array}
                  \right).$$
Obviously, $\langle xy,e_1\rangle=-3$. Therefore, without
some additional condition, one can not conclude
$\langle xy,e_i\rangle\geq 0$ for all $i \in \{1,2,3\}$.

\vskip2mm We end this section with the following property of the
elements in $K$.

\begin{Proposition}\label{1022} Let $K$ be a homogeneous cone in ${\H}$ and $x\in K$ with $x=\sum_{i,j=1}^r x_{ij}.$
If $x_{kk}=0$ for some $k\in\{1,2,\cdots,r\},$ then
$$\sum_{j=1}^r x_{kj}+\sum_{i=1}^r x_{ik}=0.$$ In particular, if
$x\in K\cap \left(\bigoplus_{i\neq j}\A_{ij}\right)$, then $x=0$.
\end{Proposition}
\textbf{Proof.} Since $x\in K$, there exists $u\in \T_+$ such that
$x=u u^\ast$. Set $u=\sum_{i,j=1}^ru_{ij}$ with $u_{ij}=0$ for every
$i,j\in\{1,2,\cdots,r\}$ such that $i>j$. Direct calculation yields
$$x_{kk}=\sum_{j=k}^ru_{kj}u_{kj}^\ast.$$
Since $\langle e_k, u_{kj}u_{kj}^\ast\rangle=\langle u_{kj},
u_{kj}\rangle=\|u_{kj}\|^2$, by the assumption $x_{kk}=0$, we
obtain $$0=\langle e_k,x_{kk}\rangle=\sum_{j=k}^r\|u_{kj}\|^2.$$
This along with $u\in\T_+$ leads to $u_{kl}=0$ for every
$l\in\{1,2,\cdots,r\}.$ Therefore, we obtain
$$x_{ik}=\sum_{l=1}^r u_{il}u_{kl}^\ast=0, ~~x_{kj}=\sum_{l=1}^r u_{kl}u_{jl}^\ast=0,~~\forall
i,j\in\{1,2,\cdots,r\},$$ as desired.{\qed}

\section{$P$ and $R_0$ properties}
We first give the definitions of various $P(P_0)$ and $R_0$
properties.
\begin{Definition}\label{dp} For a continuous function $F:\H\rightarrow \H$, we say that
it has

(i) the \emph{order-$P$ property} if for any pair $x,y\in{\H}$,
$$ (x-y)\wedge_K(F(x)-F(y))\in-(K\cap K^\ast) {\rm~ and ~} (x-y)\vee_K(F(x)-F(y))\in(K+K^\ast)~\Rightarrow~ x=y;$$

(ii) the \emph{order-$P_0$ property} if  $F(x)+\varepsilon B(x)$ has
the order-$P$ property  for any $\varepsilon>0$ where
$B:\H\rightarrow \H$ is a given linear function and satisfies
$\langle x, B(x)\rangle>0, \langle xB(x),e_i\rangle\geq 0, ~\forall
x\neq 0,~\forall i\in\{1,2,\cdots,r\}$;

(iii) the \emph{$P$ property} if for any pair $x,y\in{\H}$,
\begin{eqnarray}\sum_{l\geq
j}\left([(x-y)(F(x)-F(y))]_{lj}+[(x-y)(F(x)-F(y))]_{lj}^\ast\right)-\sum_{i=1}^r[(x-y)(F(x)-F(y))]_{ii}&\in&-(K+K^\ast)~~
\nonumber\\ \Rightarrow~ x=y;& &\nonumber\end{eqnarray}

(iv) the \emph{$P_0$ property} if  $F(x)+\varepsilon B(x)$ has the
$P$ property  for any $\varepsilon>0$;

(v) the \emph{trace-$P$ property} if for any pair $x,y\in{\H}$ with
$x\neq y$,
$$\max_i \langle (x-y)(F(x)-F(y)), e_i\rangle > 0;$$

(vi) the \emph{trace-$P_0$ property} if  $F(x)+\varepsilon B(x)$ has
the trace $P$ property  for any $\varepsilon>0$;

(vii) the \emph{uniform-trace-$P$ property} if there is an
$\alpha>0$ such that for any pair $x,y\in{\H}$,
$$\max_i \langle (x-y)(F(x)-F(y)), e_i\rangle \geq\alpha \|x-y\|.$$
 In general, we may choose the above $B$ as the identity
transformation.
\end{Definition}

\vskip1.5mm\noindent\textbf{Remark 3.1} By Proposition \ref{t101}
(b),  the implication condition of the {order-$P$ property} is
equivalent to the following: for any pair $x,y\in{\H}$,
$$ (x-y)\wedge_K(F(x)-F(y))\in-(K\cap K^\ast) {\rm~ and ~} (F(x)-F(y))\vee_{K^\ast}(x-y)\in(K+K^\ast)~\Rightarrow~ x=y.$$
Note that when $K$ is self-dual, $K\cap K^\ast=K$ and $K+K^\ast=K$.
It is easy to see that all the above order-$P$ and $P$ properties
become the order-$P$ and Jordan $P$ properties given by Tao and
Gowda \cite{TG05} in the setting of SCCP, respectively. In
particular, when $\H=\R^n$ and $K=\R^n_+$, they are all the same as
the $P$ function (see Introduction).

 \vskip1.5mm\noindent\textbf{Remark 3.2} Using the
related definitions, we can easily verify the following one-way
implications of the properties for nonlinear transformation $F$:
\begin{eqnarray}\textrm{Strong monotonicity} \Rightarrow
\textrm{~~uniform-trace-$P$~~}&\Rightarrow &\textrm{trace-$P$}
\Rightarrow \textrm{trace-$P_0$},\nonumber\\
\textrm{Strong monotonicity}\Rightarrow \textrm{strict
monotonicity}&\Rightarrow &\textrm{trace-$P$} \Rightarrow
\textrm{$P$}
,\nonumber\\
\textrm{Monotonicity}&\Rightarrow
&\textrm{trace-$P_0$}\Rightarrow\textrm{$P_0$}.\nonumber\end{eqnarray}
Here, we say that $F$ is \textit{monotone} if $\langle
x-y,F(x)-F(y)\rangle\geq0,~~\forall~x,y \in \H;$ $F$ is
\textit{strictly~monotone} if $\langle
x-y,F(x)-F(y)\rangle>0,~~\forall~x\neq y, x,y \in \H;$ and $F$ is
\textit{strongly monotone with modulus $\mu>0$} if $\langle
x-y,F(x)-F(y)\rangle\geq\mu\|x-y\|^2,~~\forall~x,y\in\H.$

\vskip1.5mm\noindent\textbf{Remark 3.3} Observe that there are very
many possible generalizations of the definition of the order-$P$
property from symmetric cones to homogeneous cones. For instance,
for any pair of sets $\hat{K}$ and $\check{K}$  such that $K\cap
K^\ast\subseteq\hat{K},\check{K}\subseteq K+K^\ast$, we say $F$ has
the \emph{order-$P$ property with respect to $\hat{K}$ and
$\check{K}$} if
$$ (x-y)\wedge_K(F(x)-F(y))\in-\hat{K} {\rm~ and ~} (x-y)\vee_{K}(F(x)-F(y))\in\check{K}~\Rightarrow~ x=y;$$
$F$ has the \emph{order-$P_0$ property with respect to $\hat{K}$ and
$\check{K}$} if $F(x)+\varepsilon B(x)$ has the order-$P$ property
for any $\varepsilon>0$. In order to
establish the existence result of a solution to HCCP (see Theorem
\ref{t104}),  we set $\check{K}:=K+K^\ast$ from the equation
(\ref{z01}) in the proof of Lemma \ref{t106}.

\vskip1.0mm Moreover, we may define the \emph{order-$P$ property} of
$F$ by the implication
$$(x-y)\wedge_{K^\ast}(F(x)-F(y))\in-(K\cap K^\ast) {\rm~ and ~} (x-y)\vee_{K^\ast}(F(x)-F(y))\in(K+ K^\ast)~\Rightarrow~ x=y;$$
and the corresponding order-$P_0$ property. However, in this case,
we cannot guarantee  the existence result of a solution to HCCP (see
Theorem \ref{t104}).

\vskip1.5mm\noindent\textbf{Remark 3.4} Similarly, there are very
many possible generalizations of the definition of the $P$ property
to homogeneous cones. For instance, for every $\tilde{K}$ such that
$K\cap K^\ast\subseteq\tilde{K}\subseteq K+K^\ast$, we say $F$ has
the \emph{$P$ property with respect to $\tilde{K}$} if
$$\sum_{l\geq
j}\left([(x-y)(F(x)-F(y))]_{lj}+[(x-y)(F(x)-F(y))]_{lj}^\ast\right)
-\sum_{i=1}^r[(x-y)(F(x)-F(y))]_{ii}\in-\tilde{K}\Rightarrow x=y;$$
and $F$ has the \emph{$P_0$ property with respect to $\tilde{K}$} if
$F(x)+\varepsilon B(x)$ has the $P$ property  for any
$\varepsilon>0$. Clearly, from the proofs of Lemma \ref{t106} and
Theorem \ref{t104}, we obtain that for any $\tilde{K}$ if $F$ has
the {$P_0$ property with respect to $\tilde{K}$} and $R_0$
properties then the associated HCCPs have solutions. If the answer
to our question in Remark 2.2 is ``yes," with the choice of
$\tilde{K}=K+K^\ast$,
then there is some hope for proving $P\Rightarrow
P_0$.

\vskip1.5mm\noindent\textbf{Remark 3.5} In the case of SCCP, the
trace-$P$ property implies the order-$P$ property, and the order-$P$
property implies the $P$ (Jordan $P$) property. However, it is not
clear whether this is valid for HCCP.

\begin{Definition} \label{d01} A continuous function $F:\H\rightarrow \H$ is
said to have the $R_0$ property if the following condition holds:
for every sequence $\{x^{(k)}\}\subset\H$ with
$$\|x^{(k)}\|\rightarrow\infty, ~~\liminf_{k\rightarrow\infty}\frac{x^{(k)}}{\|x^{(k)}\|}\in
K,~~\liminf_{k\rightarrow\infty}\frac{F(x^{(k)})}{\|x^{(k)}\|}\in
K^\ast,$$we have
$\liminf_{k\rightarrow\infty}\frac{\max_{i}\langle
x^{(k)}F(x^{(k)}), e_i\rangle}{\|x^{(k)}\|^2}>0$.
\end{Definition}

The above definition is motivated by Definition 3.2 of Tao and Gowda
\cite{TG05}, which was originally introduced for NCP by Chen and
Harker \cite{CH}. In the setting of $\R^n$, it becomes that a
continuous function $f:\R^n\rightarrow \R^n$ has the $R_0$ property
if for every sequence $\{x^{(k)}\}\subset\R^n$ with
$$\|x^{(k)}\|\rightarrow\infty,
~~\liminf_{k\rightarrow\infty}\frac{\min_{i}x_i^{(k)}}{\|x^{(k)}\|}\geq 0,
~~\liminf_{k\rightarrow\infty}\frac{\min_i
f_i(x^{(k)})}{\|x^{(k)}\|}\geq 0,$$we have
$\liminf_{k\rightarrow\infty}\frac{\max_{i}
x_i^{(k)}f_i(x^{(k)})}{\|x^{(k)}\|^2}>0$.  Clearly, when $f$ is
linear the above condition becomes equivalent to the statement that
the standard linear complementarity problem $LCP(f,0)$ has a unique
solution, namely, zero.

Applying the related definitions we can easily derive two conditions
under which the $R_0$ property holds, which is a generalization of
Proposition 3.2 in \cite{TG05}.

\begin{Proposition}\label{t105-R0} Let $F:\H\rightarrow \H$ be a continuous function. If $F$ has either the
uniform-trace-$P$ property or satisfies the following implication:
for every sequence $\{x^{(k)}\}\subset\H$ with
$$\|x^{(k)}\|\rightarrow\infty, ~~\liminf_{k\rightarrow\infty}\frac{x^{(k)}}{\|x^{(k)}\|}\in
K,~~\liminf_{k\rightarrow\infty}\frac{F(x^{(k)})}{\|x^{(k)}\|}\in
K^\ast,$$we have $\liminf_{k\rightarrow\infty}\frac{\langle x^{(k)},
F(x^{(k)})\rangle}{\|x^{(k)}\|^2}>0$,
 then  $F$ has the $R_0$ property.
\end{Proposition}

It is well-known that the notion of $R_0$ property of a function is
closely related
 its coercivity, which plays a central role in describing the boundedness of  the solution set to NCP, see, e.g.
 \cite{FP03}. In the case of HCCP, we have a similar result.

\begin{Proposition}\label{t105} Let $F:\H\rightarrow \H$ be a continuous function. If $F$ has
the $R_0$ property, then for every $\delta>0$, the set $\{x\in\H:
~x ~{\rm solves}~ HCCP(F,q),~\|q\|\leq\delta\}$ is bounded.
\end{Proposition}
\textbf{Proof.} Suppose the set $\{x\in\H: ~x ~{\rm solves}~
HCCP(F,q),~\|q\|\leq\delta\}$ is unbounded. Then, there exist
sequences $\{q^{(k)}\}$ with $\|q^{(k)}\|\leq\delta$ and
$\{x^{(k)}\}$ with $\|x^{(k)}\|\rightarrow \infty$ such that
$$x^{(k)} \in K, ~y^{(k)}=F(x^{(k)})+q^{(k)} \in {K^\ast}
,~\langle x^{(k)},y^{(k)}\rangle =0,~\forall k.$$ Since
$\{x^{(k)}\} \subset K,~\{y^{(k)}\} \subset K^\ast$ and $K$,
$K^\ast$ are closed,
$\liminf_{k\rightarrow\infty}\frac{x^{(k)}}{\|x^{(k)}\|}\in K$ and
$\liminf_{k\rightarrow\infty}\frac{y^{(k)}}{\|x^{(k)}\|}\in
K^\ast.$ Since $q^{(k)}$ is bounded,
$\liminf_{k\rightarrow\infty}\frac{q^{(k)}}{\|x^{(k)}\|}=0$. Thus,
$$\liminf_{k\rightarrow\infty}\frac{F(x^{(k)})}{\|x^{(k)}\|}=\liminf_{k\rightarrow\infty}\frac{F(x^{(k)})+
q^{(k)}}{\|x^{(k)}\|}=\liminf_{k\rightarrow\infty}\frac{y^{(k)}}{\|x^{(k)}\|}\in
K^\ast.$$ This together with the $R_0$ property of $F$ gives
$$\liminf_{k\rightarrow\infty}\frac{\max_{i}\langle
x^{(k)}F(x^{(k)}),e_i\rangle}{\|x^{(k)}\|^2}>0.$$ However, noting
that $\langle x^{(k)},y^{(k)}\rangle =0$ and the boundedness
$q^{(k)}$, by Proposition \ref{t102-2}, we obtain that for every
$i\in\{1,2,\cdots,r\}$,
$$ \frac{\langle x^{(k)}
F(x^{(k)}),e_i\rangle}{\|x^{(k)}\|^2}
=\frac{\langle
x^{(k)}y^{(k)},e_i\rangle}{\|x^{(k)}\|^2}-\frac{\langle
x^{(k)}q^{(k)},e_i\rangle}{\|x^{(k)}\|^2}=-\frac{\langle
x^{(k)}q^{(k)},e_i\rangle}{\|x^{(k)}\|^2}\rightarrow 0~{\rm as
}~k\rightarrow\infty.$$ This is a contradiction and hence the
desired conclusion follows. {\qed}

Before stating our main result in this section, we recall below a
useful result from degree theory. The topological degree technique
plays an important role in the study of complementarity problems and
variational inequality problems, see, e.g, \cite{Gow93, Ha87, HS83,
IBK97, TG05,ZH99,ZL01}. Let $\Omega$ be a bounded open set in $\H$
with its closure ${\rm cl}{(\Omega)}$ and boundary $\partial\Omega$.
For a continuous function $\Phi: {\rm cl} ({\Omega})\rightarrow \H$
and $p\not\in\Phi(\partial\Omega)$, we denote ${\rm
deg}(\Phi,\Omega,p)$ the (topological) \emph{degree} of $\Phi$ with
respect to $\Omega$ at $p$, see Lloyd \cite{Ll78} for the details.

\begin{Lemma}\label{t106-0} (Theorem 2.1.2, \cite{Ll78}) (1) Suppose that $\Phi, \varphi: {\rm
cl} ({\Omega})\rightarrow \H$ are continuous and $p\not\in
\Phi(\partial\Omega)$. If $\sup_{x\in{\rm
cl}(\Omega)}\|\Phi(x)-\varphi(x)\|< {\rm
dist}(p,\Phi(\partial\Omega))$, then  ${\rm
deg}(\varphi,\Omega,p)$ is defined and $${\rm
deg}(\varphi,\Omega,p)={\rm deg}(\Phi,\Omega,p).$$

(2) If $g_t(x)$ is a homotopy and $p\not\in g_t(\partial\Omega)$
for $0\leq t\leq 1$, then ${\rm deg}(g_t,\Omega,p)$ is independent
of $t\in [0,1]$.
\end{Lemma}

The next lemma relies heavily on the invariance of degree under
suitable homotopies and generalizes Theorem 3.1 of \cite{TG05} and
its proof.

\begin{Lemma} \label{t106}  Let $F:\H\rightarrow \H$ be a continuous function, and for every $\delta>0$ the set
\begin{eqnarray}\label{004}\{x\in\H: ~x ~{\rm solves}~ HCCP(F,q),~\|q\|\leq\delta\}\end{eqnarray} is
bounded.  If $F$ has  either the order-$P_0$ property, or the $P_0$
property, then for every $q\in \H$, the solution set of HCCP(F,q) is
nonempty and bounded.
\end{Lemma}
\textbf{Proof.} Choose any $q\in \H$. Consider the function
$$\Phi(x):=x\wedge_K (F(x)+q).$$
Define the homotopy $$G_1(x,t):=x\wedge_K
[F(x)+tq+(t-1)F(0)],~t\in [0,1].$$ Clearly, $G_1(x,0)=x\wedge_K
[F(x)-F(0)]$ and $G_1(x,1)=\Phi(x)$ for all $x$. By the
assumption, the sets $\{x\in\H: G_1(x,t)=0\}$ ($t\in[0,1]$) are
uniformly bounded. Thus, we may take a bounded open set
$\Omega\in\H$ such that
$$\bigcup_{t\in[0,1]}\{x\in\H: G_1(x,t)=0\}\subseteq \Omega.$$
Then, $0\in\Omega$ and $0\notin G_1(\partial \Omega,0)$ since $
G_1(0,0)=0.$ Therefore, by Lemma \ref{t106-0} (2),
$${\rm deg}(G_1(\cdot,0),\Omega,0)={\rm deg}(G_1(\cdot,1),\Omega,0)={\rm
deg}(\Phi,\Omega,0).$$ Define $\varphi_\varepsilon(x):=x\wedge_K
[F(x)+\varepsilon B(x)-F(0)]$ for any $\varepsilon>0$,
where $B$ is
linear and strictly monotone. Note that
\begin{eqnarray}
\|\varphi_\varepsilon(x)-G_1(x,0)\|&=&\left\|x\wedge_K [F(x)+\varepsilon
B(x)-F(0)]- x\wedge_K [F(x)-F(0)]\right\| \nonumber\\
&=& \left\|\Pi_K[x-(F(x)+\varepsilon B(x)-F(0))]-\Pi_K [x-
(F(x)-F(0))]\right\| \nonumber\\
&\leq& \left\|[x-(F(x)+\varepsilon B(x)-F(0))]-[x-
(F(x)-F(0))]\right\| \nonumber\\
&= &\|\varepsilon B(x)\|, \nonumber\end{eqnarray} where the
inequality follows from (\ref{100-3}). Since ${\rm dist}(0,
G_1(\partial \Omega,0))>0$ by $0\notin G_1(\partial \Omega,0)$, we
pick $\varepsilon_0>0$ such that
$$\sup_{x\in{\rm cl}(\Omega)}\|\varphi_\varepsilon(x)-G_1(x,0)\|<{\rm dist}(0, G_1(\partial
\Omega,0)).$$ Then, by Lemma \ref{t106-0} (1), ${\rm
deg}(G_1(\cdot,0),\Omega,0)={\rm
deg}(\varphi_\varepsilon,\Omega,0).$ So, we obtain
\begin{eqnarray}\label{002}{\rm deg}(\varphi_\varepsilon,\Omega,0)={\rm
deg}(\Phi,\Omega,0).\end{eqnarray} For small $\varepsilon>0$, we
define the homotopy
$$G_2(x,t):=x\wedge_K
[t(F(x)-F(0)+\varepsilon B(x))+(1-t)x],~t\in [0,1].$$Clearly,
$G_2(x,0)=x\wedge_K x=x$ and $G_2(x,1)=\varphi_\varepsilon(x)$ for
all $x$. We now show that $0\notin G_2(\partial \Omega,t)$ for any
$t\in[0,1]$. Suppose not, then there exist $t_0\in[0,1]$ and
$x_0\in\partial \Omega$ such that $G_2(x_0,t_0)=0$. If $t_0=0$,
then $G_2(x_0,0)=0$ means that $x_0=0$, which contradicts
$0\in\Omega$. We may assume $t_0\in (0,1]$. By Proposition
\ref{t102-2}, $G_2(x_0,t_0)=0$ is equivalent to the following
$$x_0\in K,~F(x_0)-F(0)+\varepsilon
B(x_0)+\frac{1-t_0}{t_0}x_0\in K^\ast,~\left\langle
x_0,F(x_0)-F(0)+\varepsilon
B(x_0)+\frac{1-t_0}{t_0}x_0\right\rangle=0.$$ Letting
$\widetilde{F}(x):=F(x)+\varepsilon B(x)+(\frac{1}{t_0}-1)x$, the
above can be written as
$$x_0\in K,~\widetilde{F}(x_0)-\widetilde{F}(0)\in K^\ast,~\left\langle
x_0,\widetilde{F}(x_0)-\widetilde{F}(0)\right\rangle=0.$$ Thus, by
Propositions \ref{t101} and \ref{t102-2}, we obtain
\begin{eqnarray} (x_0-0)\wedge_K
(\widetilde{F}(x_0)-\widetilde{F}(0))&=&0\in -(K\cap
K^\ast),\nonumber\\
\label{z01}~~(x_0-0)\vee_K
(\widetilde{F}(x_0)-\widetilde{F}(0))&=&x_0+
(\widetilde{F}(x_0)-\widetilde{F}(0))\in K+ K^\ast,\end{eqnarray}
and
\begin{eqnarray} \label{z02}~~[(x_0-0)(\widetilde{F}(x_0)-\widetilde{F}(0))]_{lj}=0,~\forall l,j\in\{1,2,\cdots,r\}~ {\rm
such ~that}~l\geq j.
\end{eqnarray}
If $F$ has the order-$P_0$ property, then $\widetilde{F}$ has the
order-$P$ property. Hence, by (\ref{z01}), $x_0=0$, a contradiction.
If $F$ has the $P_0$ property, then $\widetilde{F}$ has the $P$
property and hence, by (\ref{z02}), $x_0=0$. This is also a
contradiction.

Thus, $0\notin G_2(\partial \Omega,t)$. Again, by Lemma
\ref{t106-0} (2), \begin{eqnarray}\label{003}{\rm
deg}(G_2(\cdot,0),\Omega,0)={\rm deg}(G_2(\cdot,1),\Omega,0)={\rm
deg}(\varphi_\varepsilon,\Omega,0).\end{eqnarray} Notice that
${\rm deg}(G_2(\cdot,0),\Omega,0)=1$. This together with
(\ref{002}) and (\ref{003}) yields ${\rm deg}(\Phi,\Omega,0)=1$,
which says that $\Phi(x)=0$ has a solution. By Proposition
\ref{t102-2}, we proved that HCCP$(F,q)$ has a solution. The
desired conclusion follows from the assumption (\ref{004}). {\qed}

We state below our main result in this section.

\begin{Theorem}\label{t104} Let $F:\H\rightarrow \H$ be a continuous function. Suppose that
$F$ has either the order-$P_0$ and $R_0$ properties, or the $P_0$
and $R_0$ properties. Then for every $q\in\H$, the solution set of
HCCP(F,q) is nonempty and bounded.
\end{Theorem}
\textbf{Proof.} It follows immediately from
 Proposition \ref{t105} and Lemma \ref{t106}. {\qed}

As a direct consequence of the above theorem, we have the
following.

\begin{Corollary}\label{t107} Let $F:\H\rightarrow \H$ be a continuous function. Suppose that
 $F$ has either the order-$P_0$ and $R_0$ properties, or the $P_0$ and $R_0$ properties. Then, there exists $\bar{x}\in \H$ such that
 \begin{eqnarray} \bar{x}\in {\rm int}(K),~~~F(\bar{x})\in {\rm
 int}(K^\ast).\nonumber\end{eqnarray}
\end{Corollary}
\textbf{Proof}. It follows from Theorem \ref{t104} that for every
$q\in\H$, the solution set of HCCP$(F,q)$ is  nonempty and
bounded. Take $-q_0\in{\rm
 int}(K^\ast)$. Let $\hat{x}$ be a solution to HCCP$(F,q_0)$. Then,
 we have $\hat{x}\in K, \hat{y}:=F(\hat{x})+q_0\in K^\ast.$ Therefore,
 $F(\hat{x})=-q_0+\hat{y}\in {\rm  int}(K^\ast)$ and the desired
 conclusion follows from the continuity of $F$.  {\qed}

\section{GUS Property}

We say that a continuous transformation $F:{\H}\rightarrow{\H}$ has
the \emph{GUS property} if for every $q\in \H$,  HCCP$(F,q)$ has a
unique solution.  In this section, we show that the trace-$P$ property is
sufficient for $F$ having the GUS property in the setting of HCCP.
We also establish an error bound for the HCCP with the
uniform-trace-$P$ property.

\subsection{Sufficient conditions for the GUS property}

Next, we show that if $F$ has the trace-$P$ property together with
monotonicity or some other suitable property,
then the associated HCCP$(F,q)$ has
the GUS property.

\begin{Theorem}\label{t104-1} Let $F:\H\rightarrow \H$
be a monotone continuous function. Suppose that
 $F$ has  the trace-$P$ property. Then for any
$q\in\H$, the solution set of HCCP$(F,q)$ is a singleton,
if it is nonempty.
\end{Theorem}
\textbf{Proof.} Assume that there exist two solutions $x,y$ to
HCCP$(F,q)$. Then, by Proposition \ref{t102-2},
\begin{eqnarray}\label{z03}x,y\in K, ~F(x)+q,F(y)+q\in K^\ast,~\langle
x(F(x)+q),e_i\rangle=\langle y(F(y)+q),e_i\rangle=0,~\forall
i\in\{1,2,\cdots,r\}.\nonumber\end{eqnarray} Then,\begin{eqnarray}
\langle x-y,F(x)-F(y)\rangle= \langle
x-y,(F(x)+q)-(F(y)+q)\rangle=-\langle x,F(y)+q\rangle- \langle
y,F(x)+q\rangle\leq 0.\nonumber\end{eqnarray}Also, since $F$ is
monotone, $\langle x-y, F(x)-F(y)\rangle\geq 0.$ So, $\langle x-y,
F(x)-F(y)\rangle= 0$ and $\langle x,F(y)+q\rangle= \langle
y,F(x)+q\rangle= 0$. Similarly,  by Proposition \ref{t102-2},
$\langle x(F(y)+q),e_i\rangle=\langle y(F(x)+q),e_i\rangle=0$. Thus,
we obtain that
\begin{eqnarray} \langle (x-y)(F(x)-F(y)),e_i \rangle=\langle
(x-y)[(F(x)+q)-(F(y)+q)],e_i\rangle=0.\nonumber\end{eqnarray}Therefore,
the trace-$P$ property of $F$ yields $x=y$, as desired. {\qed}

As a direct application of the above Theorems \ref{t104-1} and
\ref{t104} and the connection between various $P$ properties, we
have that under the monotonicity assumption,
the uniform-trace-$P$ property implies the GUS property.
\begin{Corollary}\label{t104-2}
Let $F:\H\rightarrow \H$ be a monotone continuous function.
If  $F$ has  the uniform-trace-$P$ property then it has the GUS
property.
\end{Corollary}

Note that in the special
case of SCCP, if $\bar{x}$ is a  solution to SCCP$(F,q)$,
then $\bar{x}, F(\bar{x})+q$ share a common Jordan frame.
This  motivates us to give the following result.

\begin{Theorem}\label{t104-1-1} Let $F:\H\rightarrow \H$ be a continuous function. Suppose that
 $F$ has  the trace-$P$ property. For every $q\in\H$, if  HCCP$(F,q)$ has
 a solution $x$ such that $x=\sum_{i=1}x_ie_i, F(x)+q=\sum_{i=1}y_ie_i$, then
it has a unique solution.
\end{Theorem}
\textbf{Proof.} As in the proof of Theorem \ref{t104-1}, assume that there exist two solutions $x,y$ to
HCCP$(F,q)$. Then, by Proposition \ref{t102-2},
\begin{eqnarray}\label{z03}x,y\in K, ~F(x)+q,F(y)+q\in K^\ast,~\langle
x(F(x)+q),e_i\rangle=\langle y(F(y)+q),e_i\rangle=0,~\forall
i\in\{1,2,\cdots,r\}.\nonumber\end{eqnarray} Thus, we obtain that
\begin{eqnarray} \langle (x-y)(F(x)-F(y)),e_i \rangle&=& \langle
(x-y)[(F(x)+q)-(F(y)+q)],e_i\rangle\nonumber\\
&=&-\langle x(F(y)+q),e_i\rangle- \langle
y(F(x)+q),e_i\rangle\nonumber\\
&\leq & 0,\nonumber\end{eqnarray}where the second equality holds by
Proposition \ref{t102-2} and the inequality holds by Proposition
\ref{t103}. Therefore, the trace-$P$ property of $F$ yields $x=y$, as
desired. {\qed}

\begin{Corollary}\label{t104-2-1} Let $F:\H\rightarrow \H$ be a
continuous function.
If $F$ has  the uniform-trace-$P$ property and for any $q\in\H$,
HCCP$(F,q)$ has a solution
$x$ such that $x=\sum_{i=1}x_ie_i, F(x)+q=\sum_{i=1}y_ie_i$,
then $x$ is the unique solution to HCCP$(F,q)$.
\end{Corollary}

\subsection{Error Bound}

\vskip1mm We give an error bound for the HCCP with the
uniform-trace-$P$ property.

\begin{Theorem} Suppose $F$ has the uniform-trace-$P$ property with modulus $\alpha>0$ and is Lipschitz continuous
with constant $\kappa>0$.  Let $x^\ast$ be the unique solution of
problem HCCP$(F,q)$ such that $x^\ast=\sum_{i=1}x_i^\ast e_i, F(x^\ast)+q=\sum_{i=1}y_i^\ast e_i$. Then
\begin{eqnarray}\label{706}\frac1{2+\kappa}\|x\wedge_K (F(x)+q)\|\leq \|x-x^\ast\|\leq
\frac{1+\kappa}{\alpha}\|x\wedge_K (F(x)+q)\|,~\forall x\in
\H.\end{eqnarray}
\end{Theorem}
\textbf{Proof.} Notice that $F(x)+q-x\wedge_K
(F(x)+q)=\Pi_{K^\ast}(F(x)+q-x)\in K^\ast,$~ $x-x\wedge_K
(F(x)+q)=\Pi_K(x-F(x)-q)\in K,$  and $x^\ast\wedge_K
(F(x^\ast)+q)=0$. By Propositions \ref{t102-2} and \ref{t103} and
the fact that $x^\ast$ solves HCCP$(F,q)$ such that $x^\ast=\sum_{i=1}x_i^\ast e_i, F(x^\ast)+q=\sum_{i=1}y_i^\ast e_i$, we obtain
\begin{eqnarray} \langle[F(x)+q-x\wedge_K (F(x)+q)][x^\ast-x^\ast\wedge_K
(F(x^\ast)+q)], e_i\rangle&\geq& 0,\nonumber\\ \langle[x-x\wedge_K
(F(x)+q)][F(x^\ast)+q-x^\ast\wedge_K (F(x^\ast)+q)],
e_i\rangle&\geq& 0,\nonumber\\ \langle [F(x^\ast)+q-x^\ast\wedge_K
(F(x^\ast)+q)][x^\ast-x^\ast\wedge_K
(F(x^\ast)+q))],e_i\rangle&=&0.\nonumber\end{eqnarray} Thus, by
direction calculation, we obtain that for every
$i\in\{1,2,\cdots,r\}$,
\begin{eqnarray}
& &\langle[F(x)+q-x\wedge_K (F(x)+q)][x-x\wedge_K (F(x)+q)], e_i\rangle\nonumber\\[1.5mm]
&\geq&\langle([F(x)+q-x\wedge_K
(F(x)+q)]-[F(x^\ast)+q-x^\ast\wedge_K
(F(x^\ast)+q)])\nonumber\\[1.5mm]& &
([x-x\wedge_K (F(x)+q)]-[x^\ast-x^\ast\wedge_K (F(x^\ast)+q)],e_i\rangle\nonumber\\[1.5mm]
&=&\langle [F(x)-F(x^\ast)-x\wedge_K (F(x)+q)][x-x^\ast-x\wedge_K (F(x)+q)],e_i\rangle\nonumber\\[1.5mm]
&\geq &\langle [F(x)-F(x^\ast)](x-x^{\ast}),e_i\rangle -\langle
x\wedge_K (F(x)+q)[F(x)-F(x^\ast)+x-x^\ast],e_i\rangle\nonumber\\[1.5mm]
&\geq &\langle (F(x)-F(x^\ast))(x-x^\ast),e_i\rangle- \|x\wedge_K
(F(x)+q)\|\cdot (1+\kappa)\|x-x^\ast\|,\nonumber\end{eqnarray} where
the second inequality holds by the fact $\langle
a^2,e_i\rangle=\langle aa^{\ast},e_i\rangle\geq 0$ for all $a\in \H$
because $a=a^{\ast}$; the last inequality follows from the Lipschitz
continuity of $F$ and $\|ab\|\leq\|a\|\|b\|$ for any $a,b\in\A$ by
Theorem 2 of \cite{Chua08}. By Propositions \ref{t101} and
\ref{t102-2}, $\langle [F(x)+q-x\wedge_K (F(x)+q)][x-x\wedge_K
(F(x)+q)],e_i\rangle=0$. Thus, we conclude from the above inequality
and the uniform-trace-$P$ property of $F$,
\begin{eqnarray}(1+\kappa)\|x\wedge_K (F(x)+q)\|
\|x-x^\ast\|&\geq&\max_{i}\langle
(F(x)-F(x^\ast))(x-x^\ast),e_i\rangle\geq\alpha\|x-x^\ast\|^2.\nonumber\end{eqnarray}
This leads to the right-hand side of inequality (\ref{706}).

\vskip1mm Note that $\|\Pi_K(y)-\Pi_K(z)\|\leq\|y-z\|$ for every
$y,z\in\H$. From the Lipschitz continuity of $F$, we obtain by
direct manipulation that
\begin{eqnarray}\|x\wedge_K (F(x)+q)\|&=&
\|[x-\Pi_K(x-F(x)-q)]-[x^\ast-\Pi_K(x^\ast-F(x^\ast)-q)]\|\nonumber\\[1mm]
&=& \|[x-x^\ast]-[\Pi_K(x-F(x)-q)-\Pi_K(x^\ast- F(x^{\ast})-q)]\|\nonumber\\[1mm]
&\leq& \|x-x^\ast\|+\|\Pi_K(x- F(x)-q)-\Pi_K(x^\ast-F(x^\ast)-q)\|\nonumber\\[1mm]
&\leq& \|x-x^\ast\|+\|(x-F(x)-q)-(x^\ast-F(x^\ast)-q)\|\nonumber\\[1mm]
&\leq& 2\|x-x^\ast\|+\|F(x)-F(x^\ast)\|\nonumber\\[1mm]
&\leq& (2+\kappa)\|x-x^\ast\|.\nonumber\end{eqnarray} This means
that the left-hand side of inequality (\ref{706}) holds.   {\qed}
Note that in the SCCP context, the condition  $x^\ast=\sum_{i=1}x_i^\ast e_i, F(x^\ast)+q=\sum_{i=1}y_i^\ast e_i$
is  naturally satisfied. Then, the above error bound result
becomes a generalization of the corresponding one for NCP presented by Chen and Harker
\cite{CH}.  It is not clear
whether the above error bound would
hold for every monotone, continuous function $F$ with
the uniform-trace-$P$ property.

\section{Final Remarks}

In this paper, by employing the $T$-algebraic characterization of
homogeneous cones we prove that if a continuous function has either
the order-$P_0$ and $R_0$, or the $P_0$ and $R_0$ properties then
all the associated HCCPs have solutions.
We give sufficient conditions under which
the associated HCCP has the GUS property.
Moreover, we establish a global error bound for
HCCP under some conditions.

\vskip2mm Many of our results apply to the more general setting of
arbitrary convex cones (the order-$P$ property, order-$P_0$
property, trace-$P$ property, and $R_0$ property). Further
generalizations of similar results and the theory to arbitrary
convex cone setting represent
a good direction for future research. The
design of algorithms for HCCP (and beyond) and a study of their
mathematical and computational properties provide other interesting
future research avenues.

\vskip3mm\vskip1.5mm \noindent{\Large\textbf{Acknowledgments}} 
The authors thank Nan Lu
of Tianjing University for very helpful discussions, which helped us
eliminate a very important error in the first version of this paper.

 
\end{document}